\documentclass[11pt]{article}
\usepackage{amssymb,latexsym,epsf,xypic}
\usepackage[mathscr]{euscript}

\makeatletter \@addtoreset{equation}{section} \makeatother

\newtheorem{Lemma}[equation]{Lemma}
\newtheorem{Theorem}[equation]{Theorem}

\newenvironment{Proof}{\noindent\emph{Proof\ }}{\hfill$\square$\\}

\newcommand\C{\mathbb C\,}
\newcommand\Z{\mathbb Z}
\newcommand\Pee{\mathbb P}

\newcommand\OO{\mathscr O}
\newcommand\I{\mathcal I}
\newcommand\To{\longrightarrow}

\newcommand\into{\hookrightarrow}
\newcommand\res{\arrowvert_}

\newcommand\id{\mathrm{id}}
\newcommand\End{\mathrm{End}\,}

\begin{document}

\title{Examples of bundles on Calabi-Yau 3-folds for string theory
compactifications}
\author{R.\,P. Thomas}
\date{}
\maketitle

Witten asked in 1997 (and no doubt much earlier too) if there was an
example of a rank four bundle $E$ on a Calabi-Yau 3-fold $X$
satisfying the following conditions:
\begin{itemize}
\item $X$ is not simply-connected,
\item $c_1(E)=0$; in fact $\Lambda^4E\cong\OO_X$,
\item the holomorphic Euler characteristic $\chi(E)$ of $E$ equals 3,
\item $E$ is slope-polystable (a direct sum of slope-stable bundles of
degree zero with respect to the K\"ahler form), i.e. admits an $SU(4)$
Hermitian-Yang-Mills connection by the Donaldson-Uhlenbeck-Yau theorem,
\item the pairing $H^1(\Lambda^2E)\otimes H^1(E)\otimes H^1(E)\to\C$ (given
by cup-product, wedging, and $H^3(\Lambda^4E)\cong H^3(\OO_X)\cong\C$)
should be non-zero,
\item the same pairing $H^1(\Lambda^2E\otimes\alpha)\otimes
H^1(E\otimes\beta)\otimes H^1(E\otimes\gamma)\to\C$, where
$\alpha,\,\beta,\,\gamma$ are representations of the fundamental group
of $X$ such that $\alpha\otimes\beta\otimes\gamma\cong\OO_X$,
vanishes for $\alpha,\,\beta$ and $\gamma$ non-trivial,
\item the last two conditions also hold with $E$ replaced by $E^*$
throughout,
\item one can analyse the pairings $H^1(\End\,E)\otimes H^1(E)\otimes
H^1(E^*)\to H^3(\OO_X)\cong\C$, and
\item $c_2(X)-c_2(E)$ is $c_2(F)$ for some slope-polystable bundle
$F$ of any small rank and $c_1(F)=0$.
\end{itemize}

This has something to do with compactifying a 10-dimensional string theory
on $X$ to the supersymmetric standard model in four dimensions to study
the half-life of the proton. $E$ and $F$ are to be embedded in $E_8$
bundles, and $\chi(E)=3$ is the ``generation number'' -- the number of
families of quarks. This is the limit of my understanding, but fortunately
Witten distilled the physics down to the above purely
mathematical question.

Physicists usually concentrate on producing examples satisfying as many of
the topological conditions as possible (see \cite{DOPW} for the current
state of the art), hoping that if the moduli are big enough there will
be at least one bundle satisfying the conditions on the pairings.
(The first example satisfying these topological conditions was given
by Tian and Yau in \cite{Y}\,\cite{TY}, but their example with
$E=TX\oplus\OO$
does not satisfy the other constraints.) Physicists also tend to use the
Friedman-Morgan-Witten method of constructing bundles on elliptically
fibred 3-folds. Here the Serre construction is used; in fact even though
the 3-folds below are elliptic the FMW method \emph{does not apply}
because the bundles \emph{are not stable on the elliptic fibres} and the
nonexistence results of physicists (e.g. \cite{DLOW} rules out these
3-folds) do not apply directly. Here the problem is tackled
from the other end, constructing bundles satisfying the pairing
conditions (which I feel ought to be the most difficult) with enough freedom
to try and control the topological ones. This is only partially successful
-- the last condition still eludes me (and so \cite{DLOW} may yet apply to
this case). $c_2(E)$ is kept as low as possible
to make satisfying the last condition feasible, but I have yet to find
the required $F$. The necessary condition (the Bogomolov inequality
$c_2(F)\,.\,\Omega\ge0$, where $\Omega$ is the K\"ahler form) for the
existence of a stable $F$ is easily satisfied, but
$c_2(F)$ is \emph{not} effective. I still believe the class $c_2(F)$ may be
represented by a stable $F$ but my attempts to find one over the last 2
years have failed. (Since it is usually assumed in the
physics literature that $c_2$ of a stable bundle, with $c_1=0$, must be
an effective curve we give an example to show that this need not be true.)

The first example below is on $K3\times T^2$; when I showed this to
Witten another condition was promptly added to the above list -- that
$X$ should not be $K3\times T^2$ (in fact that the holonomy of $X$ should
be bigger than $SU(2)\subset SU(3)$). It is included below and worked out
in full, however, as it displays most of the ideas of, and is good
motivation for, the other two examples, which we run through more briefly
as most of the principles are the same. These second two examples take
place on the $SU(2)\times\Z_2$ holonomy manifold $(K3\times T^2)/\Z_2$
where the $\Z_2$-action is an Enriques action on the $K3$ times by $-1$ on
$T^2$, and so is free and preserves the canonical class. The bundles we
find are not really full $SU(4)$-bundles, but have smaller structure
group $U(2)\times_{U(1)}U(2)\subset SU(4)$ (i.e. they are direct sums of
rank two bundles of opposite $c_1$); I think Witten's intention
was to deform them to non-split bundles -- the pairings involving the
deformation space $H^1(\End E)$ would enable one to study such deformations.
This paper is inevitably a rather dry list of mathematical constructions,
but the general technique used to satisfy the conditions on the pairings
should be clear from the first example.

\noindent \textbf{Acknowledgements.} 
Many thanks to Edward Witten for discussing the problem with me several
times and allowing me to use his results and intuition; I have probably
misquoted him and of course have not solved the problem anyway.
Thanks to Ian Graham for the use of his house in the south of
France, and to Ian Dowker and Ivan Smith for $K3$ advice.
This work was mainly done at the Institute for Advanced Study in
1997--98 under NSF grant DMS 9304580. Having finally given up on
solving the last condition, a broken man, I have decided to circulate
this in the hope that someone else might use the same techniques to
better effect. I am grateful to Professors Yau and Taubes for their
current support at Harvard University where this paper was written up.
Thanks to Tony Pantev for his encouragement and interest.

\section{$K3\times T^2$}

Let $S$ be a smooth $K3$ $(2,3)$-divisor in $\Pee^1\times\Pee^2$.
Let $\omega_1,\,\omega_2$ be (the restrictions to $S$ of) the pullbacks
to $\Pee^1\times\Pee^2$ of the Fubini-Study K\"ahler forms on $\Pee^1,\,
\Pee^2$ respectively. They are the first Chern classes of the line bundles
$\OO(1,0)$ and $\OO(0,1)$, in the obvious notation. Tensor powers of these
give the line bundles $\OO(i,j)$. Note that $\omega_1^2=0,\ \omega_1
\omega_2=3,\ \omega_2^2=2$, and that $\omega_2$ is a nondegenerate
K\"ahler form on $S$.

The line bundle $L=\OO(-1,1)$ is of degree $-1$ (with respect to
$\omega_2$), has $c_1(L)^2=-4$, and it and its dual are acyclic -- that is
$$
H^i(L)=0=H^i(L^*) \quad \forall i.
$$
This can be seen by simple exact sequences on $\Pee^1\times\Pee^2$, or from
the obvious fact that $L$ and $L^*$ have no sections, and Riemann-Roch.

We are going to define two rank 2 bundles on $S$, with determinant
$L$ and $L^*$ respectively, by the Serre construction on a surface
(\cite{GH} p 726, \cite{DK} Chapter 10), which we briefly describe now.
Just as (codimension 1) divisors correspond to (rank 1) line
bundles, codimension 2 (i.e. dimension zero, on $S$) subschemes $Z$ 
sometimes correspond to rank 2 bundles $E$ via zero sets of
sections $s\in H^0(E)$. Suppose we have such an $s$ with zero locus $Z$.
Then, just as wedging with a non-zero vector $v\in V$
in a 2-dimensional vector space $V$ gives an exact sequence
$0\to\C\to V\to\Lambda^2V\to0$, wedging with $s\in H^0(E)$ gives a
sequence of sheaves $0\to\OO_X\to E\to\Lambda^2E\to0$
which is exact away from the zeros of $s$. If $s$ vanishes only
in codimension 2 (i.e. not on a divisor) this sequence is in fact globally
exact except at the last term, where it is clearly onto only those
sections of $\Lambda^2E$ that vanish on $Z$, so giving the exact sequence
$$
0\to\OO_X\to E\to L\otimes\I_Z\to0,
$$
where $\I_Z$ is the ideal sheaf of functions vanishing on $Z$, and $L$
is the line bundle $\Lambda^2E$.

The Serre construction provides a partial converse to this construction:
given $Z$ and $L$ one tries to reconstruct an $E$ with $\Lambda^2E\cong L$
and $s\in H^0(E)$ (with zeros on $Z$) as an extension of $L\otimes\I_Z$
by $\OO_X$ as in the above sequence. Such an extension is of course
given by an element of $H^1(L^*)$ away from $Z$, which we may think of
as an $L^*$-valued one-form that is
$\bar\partial$-closed away from $Z$. Local analysis on $Z$ (\cite{DK}
Chapter 10) shows that
the extension $E$ being locally free is equivalent to $\bar\partial$
of the form being some non-zero multiple $a_z\delta_z$ of
the Dirac delta at each point $z$ of $Z$. Thus the condition for the
global existence of such a vector bundle $E$ is that the class in
$H^2(L^*)$ defined by some combination $\sum_{z\in Z}a_z\delta_z\ \ 
(a_z\ne0\ \forall z\in Z$) is $\bar\partial$ of
something, i.e. that it is zero in cohomology. (Here the
residue data $a_z$ is really an element of the line $L^*\otimes K_X^*
\res z^{\ }$, thus giving a linear functional on $H^0(L\otimes K_X)$ by
restriction to $z$. By Serre duality $H^0(L\otimes K_X)^*\cong H^2(L^*)$
this gives an element of $H^2(L^*)$ as claimed. The $a_z$s give the dual of
the determinant of the derivative $ds\in(T^*X\otimes E)\res z^{\ }$ of
the section $s\in H^0(E)$ at the zeros $z$.)

One upshot of all this which will suffice for our needs is that
\emph{if $H^2(L^*)=0$ then $E$ and $s$ exist}.

So define the dual $A^*$ by the exact sequence of sheaves
\begin{equation} \label{A*}
0\to\OO\stackrel{s_{A^*}\,}\To A^*\to\I_2(1,-1)\to0,
\end{equation}
where $\I_2(1,-1)$ denotes the ideal sheaf of functions vanishing at two
fixed points on $S$, twisted by $L^*=\OO(1,-1)$. There is no obstruction
to defining bundles in this way (with a section $s_{A^*}$ vanishing
exactly at the two points) as $H^2(L)=0$.

From the above sequence it is evident that
$$
\Lambda^2(A)=L,\quad c_2(A^*)=2=c_2(A),\quad H^0(A)=0,\quad
H^0(A^*)=\C\,.\,s_{A^*}.
$$

Similarly we define $B$ by
\begin{equation} \label{B}
0\to\OO\stackrel{s_B\,}\To B\to\I_3(1,-1)\to0,
\end{equation}
with
$$
\Lambda^2(B)=L^*,\quad c_2(B)=3,\quad H^0(B)=\C\,.\,s_{B},\quad
H^0(B^*)=0.
$$

Now let $T$ denote an elliptic curve, and set $X=S\times
T\stackrel{\pi\,}\To S$. Finally denote by $\OO_T(n)$
the $n$th power of the pull-back to $X$ of a fixed degree one line
bundle on $T$, with first Chern class $n\omega_T$. Then we can define
\begin{equation} \label{E}
A':=\pi^*A\otimes\OO_T(3), \quad B':=\pi^*B\otimes\OO_T(-3), \quad
E:=A'\oplus B'.
\end{equation}

Assume that $S$ is a generic divisor in $\Pee^1\times\Pee^2$, so that
its only line bundles are the $\OO(i,j)$s by Noether-Lefschetz theory.
Then we have

\begin{Theorem}
The rank 4 bundle $E$ defined in (\ref{E}) satisfies all but the last of
the conditions listed at the start of this paper, with respect to the
K\"ahler form $\Omega=\pi^*\omega_2+6\omega_T$.
\end{Theorem}

\begin{Proof}
Throughout we will often suppress pull-backs for clarity; thus $A$ will
often denote $\pi^*A$ and $\omega_2$ will be confused with $\pi^*\omega_2$.

Firstly $\Lambda^4E\cong\Lambda^2A\otimes\OO_T(6)\otimes\Lambda^2B\otimes
\OO_T(-6)\cong L\otimes L^*\cong\OO_X$ is trivial. By Riemann-Roch
the holomorphic Euler characteristics $\chi(A),\,\chi(B)$
are $0$ and $-1$ respectively, so that
$$
\chi(E)=\chi(\OO_T(3))\chi(A)+\chi(\OO_T(-3))\chi(B)=3(\chi(A)-\chi(B))=3,
$$
as required. The choice of K\"ahler form $\Omega=\omega_2+6\omega_T$
ensures that $A'$ and $B'$ have degree zero, so $E=A'\oplus B'$ is
polystable if and only if $A'$ and $B'$ are stable, which in turn is
equivalent to $A$ and $B$ being stable on the $K3$ surface $S$. This
follows from

\begin{Lemma}
Recall that $S$ was chosen such that its only line bundles are the
$\OO(i,j)$s. Then letting $P$ denote either of $A^*$ or $B$ on $S$,
$P$ is stable with respect to $\omega_2$.
\end{Lemma}

\begin{Proof}
We must show that $P(i,j)$ has no sections
for $0\ge$\,deg\,$P(i,j)=6i+4j+1$, i.e. for $3i+2j\le-1$.
But the presentations (\ref{A*},\,\ref{B}) give the sequence
$$
0\to\OO(i,j)\to P(i,j)\to\I_n(i+1,j-1)\to0,
$$
where $n$ is either 2 or 3 points on $S$. The degree of $\OO(i,j)$
is $3i+2j\le-1$, so it has no sections. Similarly since the degree of
$\OO(i+1,j-1)$ is $3i+2j+1\le0$, this line bundle can only have a section
if it is trivial and the section has no zeros. Thus $\I_n(i+1,j-1)$ has
no sections for $n>0$, and we have shown that $P(i,j)$ has no sections
and so is stable.
\end{Proof}

We now turn to the pairings $H^1(\Lambda^2E\otimes\alpha)\otimes H^1(E
\otimes\beta)\otimes H^1(E\otimes\gamma)\to\C$, for $\alpha,\,\beta,\,
\gamma$ flat line bundles with $\alpha\otimes\beta\otimes\gamma\cong
\OO_X$. This is what motivated the construction
of $E$; the basic idea being that the cohomology of flat line
bundles on the elliptic curve $T$ behaves in a way that is reminiscent
of Witten's condition on the pairings; namely it is non-trivial if
and only if the line bundle is trivial.

\begin{Lemma}
Given representations $\alpha,\,\beta,\,\gamma$ of the fundamental
group of $X$, such that $\alpha\otimes\beta\otimes\gamma\cong\OO_X$,
consider the pairing $H^1(\Lambda^2E\otimes\alpha)\otimes
H^1(E\otimes\beta)\otimes H^1(E\otimes\gamma)\to\C$. Then this
is non-zero for the representations trivial, and vanishes
for $\alpha,\,\beta$ and $\gamma$ non-trivial.

The same is also true with $E$ replaced by $E^*$ throughout.
\end{Lemma}

\begin{Proof}
$H^1(\Lambda^2E\otimes\alpha)=H^1(\Lambda^2A'\otimes\alpha)\oplus
H^1(A'\otimes B'\otimes\alpha)\oplus H^1(\Lambda^2B'\otimes\alpha)$, and
the first and last terms vanish by the K\"unneth formula, since $L=
\Lambda^2A$ and $L^*=\Lambda^2B$ have no cohomology on $S$ and
$\alpha,\,\beta,\,\gamma$ are pulled up from $T$. Thus 
$H^1(\Lambda^2E\otimes\alpha)=H^1_S(A\otimes B)\otimes H^0_T(\alpha)\ \oplus
\ H^0_S(A\otimes B)\otimes H^1_T(\alpha)$ and the pairing reduces to
$$ \hspace{-25mm}
\left[H^1_S(A\otimes B)\otimes H^0_T(\alpha)\right]\ \otimes\ 
\left[H^1_S(A)\otimes H^0_T(\OO_T(3)\otimes\beta)\right] \vspace{-2mm}
$$$$
\hspace{3cm} \otimes\ \left[H^0_S(B)\otimes
H^1_T(\OO_T(-3)\otimes\gamma)\right]\ \to\ \C, \\
$$
plus the same with $\beta$ and $\gamma$ exchanged.
(The pairing involving $\left[H^0_S(A\otimes B)\otimes H^1_T(\alpha)\right]$
vanishes because no $H^1_S(B)$ term survives the Leray spectral sequence
to cup it with.)

So we see that for $\alpha$ non-trivial the whole pairing vanishes. (Note
that for $\alpha$ trivial but $\beta\cong\gamma^{-1}$ non-trivial the
pairing does not vanish; I am not sure if this is relevant for the
physics.) For $\alpha=\beta=\gamma=\OO$ we are left with showing, then,
that the pairing \emph{on $S$}
$$
H^1_S(A\otimes B)\otimes H^1_S(A)\otimes H^0_S(B)\to\C
$$
is non-zero. (The full pairing on $E$ is two copies of the tensor product
of this with the non-vanishing cup-product $H^0_T(\OO_T(3))\otimes
H^1_T(\OO_T(-3))\to\C$.)

But this pairing is Serre-dual to $H^1(A)\otimes H^0(B)\to H^1(A\otimes
B)$, and $H^0(B)=\C\,.\,s_B$, so it is sufficient to show that the map
$$
H^1(A)\stackrel{s_B\,}\To H^1(A\otimes B)
$$
is non-zero.

Tensoring (\ref{B}) with $A$ (recalling that $A\otimes\Lambda^2A^*\cong
A^*$) and taking cohomology gives
$$
H^0(A^*\otimes\I_3)\to H^1(A)\stackrel{s_B\,}\To H^1(A\otimes B).
$$
(\ref{A*}) shows that $H^0(A^*\otimes\I_3)=0$ and $H^1(A)\cong
H^1(\I_2)\cong\,\mathrm{coker}\,[H^0(\OO_S)\to H^0(\OO_2)]\cong\C$,
so the second map in the above sequence does not vanish, as required.

The dual pairing, with $E$ replaced by $E^*$, is similar. As above,
since $L$ and $L^*$ have no cohomology on $S$, and $\alpha$ has no cohomology
on $T$ unless it is trivial, the pairing vanishes for $\alpha$ non-trivial,
and in the $\alpha=\beta=\gamma=\OO$ case it quickly reduces to
\begin{equation} \label{pair}
H^1_S(A^*\otimes B^*)\otimes H^0_S(A^*)\otimes H^1_S(B^*)\to\C,
\end{equation}
with $H^0_S(A^*)$ generated by $s_{A^*}$. Now (\ref{A*})
twisted by $B^*$ yields
$$
H^0(B\otimes\I_2)\to H^1(B^*)\stackrel{s_{A^*}\,}\To H^1(A^*\otimes B^*).
$$
Using (\ref{B}) we can see that the first group is either $\C$ or 0
(depending on whether or not the 2 points used to define $A^*$ are
a subset of the 3 points used to define $B$) and the second group is $\C^2$.
Thus the second map, which is Serre-dual to (\ref{pair}), has rank
1 or 2, and the dual pairing is non-zero also.
\end{Proof}

Finally then we want to understand $H^1(\End\,E)\otimes H^1(E)\otimes
H^1(E^*)\to\C$. Using the fact that $H^0_T(\OO_T(6))\cong\C^6$ etc.,
we can express $H^1(\End E)$ in terms of cohomology groups on $S$ as
\begin{equation} \label{decomp}
H^1(\End A)\ \oplus\ H^1(\End B)\ \oplus\ H^0(A^*\otimes B)\otimes\C^6
\ \oplus\ H^1(B^*\otimes A)\otimes\C^6.
\end{equation}
This is easily computed to be $6+10+1.6+12.6=94$-dimensional. But we
will find that only the first two groups (16 dimensions) contribute to
the pairing. The part of the pairing involving the last group in the
above decomposition is
$$
H^1((B')^*\otimes A')\otimes H^1(B')\otimes H^1((A')^*)\to\C,
$$
which, by the K\"unneth formula, reduces to a pairing on $S$ tensored
with the cup-product
$$
H^0_T(\OO_T(6))\otimes H^1_T(\OO_T(-3))\otimes H^1_T(\OO_T(-3))
$$
on $T$, and this vanishes.

The third group in (\ref{decomp}) is involved in the cup-product pairing
$$
H^0(A^*\otimes B)\otimes H^1(A)\otimes H^1(B^*)\to\C
$$
on $S$ (tensored with a pairing on $T^2$). We will now show that this
vanishes.

Tensoring (\ref{B}) with $A^*$ shows that $H^0(A^*\otimes B)$ is spanned
by $s_{A^*}\otimes s_B$, so it is enough, by Serre-duality, to show the
vanishing of
$$
\diagram
H^1(B^*)\rto^{s_{A^*}\otimes s_B} & H^1(A^*).
\enddiagram
$$
(\ref{B}) twisted by $L$ shows that wedging with $s_B$ gives an
isomorphism $H^1(B^*)\stackrel{\wedge s_B}\To H^1(\I_3)$,
while (\ref{A*}) shows that the quotient of $\I_3\stackrel{s_{A^*}}\To
A^*$ is $\I_2(1,-1)\oplus\OO_3$, giving an exact sequence
$$
0\to\I_3\stackrel{s_{A^*}}\To A^*\to\I_2(1,-1)\oplus\OO_3\to0,
$$
whose cohomology gives the exact sequence
$$
\diagram
0\To\C\To\C^3 \rto & H^1(\I_3) \rto^{s_{A^*}} & H^1(A^*) \rto
& \C^2\To\C\To0. \\
& H^1(B^*) \uto_{\wedge s_B}^\wr\urto_{s_{A^*}\otimes s_B}
\enddiagram
$$
It is easy to check that the diagram commutes, and $H^0(\I_3)\cong\C^3$,
so the map labelled $s_{A^*}$ must be zero; therefore that marked
$s_{A^*}\otimes s_B$ also vanishes, as claimed.

So from the decomposition (\ref{decomp}) we see that the coupling
$H^1(\End\,E)\otimes H^1(E)\otimes H^1(E^*)\to\C$ on $X$ reduces to
the direct sum of the corresponding couplings for $A'$ and $B'$. By
the K\"unneth formula and Serre duality we are left with understanding
the couplings
$$
\begin{array}{ccc}
H^1(A)\otimes H^0(A^*)\to H^1(\End A) & \quad\mathrm{and}\quad &
H^0(B)\otimes H^1(B^*)\to H^1(\End B) \vspace{2mm} \\
\quad\C\otimes\C\ \to\ \C^6 && \quad\,\C\otimes\C^2\ \to\ \C^{10}.
\end{array}
$$
Twisting (\ref{A*}) by $A$ and taking cohomology gives an exact sequence
\begin{equation} \label{p}
\ldots\to H^1(A)\stackrel{s_{A^*}\,}\To H^1(\End A)\to
H^1(A^*\otimes\I_2)\to\ldots
\end{equation}
Here the first map is injective because the previous two terms in the
sequence are $H^0(\End A)\to H^0(A^*\otimes\I_2)$, which is $\C\,.\,\id
\stackrel{s_{A^*}\,}\To\C\,.\,s_{A^*}$ and so an isomorphism.
($A$ is stable so has only scalar endomorphisms; similarly the last map
can be seen to be onto and the last group
isomorphic to $\C^5$.) But this first map is the $A$-pairing, which is
therefore an injection $\C\otimes\C\into\C^6$.

As for $B$, twist (\ref{B}) by $B^*$ to get, by similar arguments,
$$
0\to H^1(B^*)\stackrel{s_B}\To H^1(\End B)\to H^1(B\otimes\I_3)\to0,
$$
so that again the pairing is an injection $\C\otimes\C^2\into\C^{10}$.
\end{Proof}

So all that is left is to find an $F$ with trivial determinant and
$c_2(F)=c_2(X)-c_2(E)=15[T]+6\omega_T(\omega_1-\omega_2)$ (where $[T]$
denotes the class in $H^4(X;\Z)$ Poincar\'e-dual to any torus fibre $T$).
Suppose we try for a rank 2 $F$, then one can calculate that
$c_2(F':=F\otimes L^*\otimes\OO_T(3))=11[T]$, which is effective.
Similarly if $F$ has rank 3 then $c_2(F':=F\otimes L^*\otimes\OO_T(1))
=4[T]$. One could therefore try to use these facts to create a stable
$F$, for instance
by using the Serre construction to manufacture a $G$ with $\Lambda^2G=L^*$
and $c_2(G)=11[T]$ then twisting by $\OO_T(n)$\,s and modifying in
codimension 1 with elementary transformations, etc, to get the desired
$F'$. All my attempts have produced \emph{unstable} bundles, however.

We note here that it is \emph{not} necessary for $c_2(F)$ to be effective
(i.e. represented by a holomorphic curve) for $F$ polystable with $c_1(F)=0$.
Indeed the polystable $E$ constructed above (\ref{E}) has $c_1(E)=0$ and
$c_2(E)=9[T]-6\omega_T(\omega_1-\omega_2)$. Thus $c_2(E)\,.\,\omega_2=-6$
is negative and $c_2(E)$ cannot be effective.

\section{$(K3\times T^2)/\Z_2$}

Now let $X$ be the $SU(2)\times\Z_2$ holonomy manifold
$(K3\times T)/\Z_2$, where the $K3$ is a universal cover of an Enriques
surface $S=K3/\sigma$ and $T$ is an elliptic curve (which therefore has a
zero and a multiplication by $-1$). Then the $\Z_2$-action is generated by
$\sigma\times(-1)$, and so is free (since $\sigma$ is) and preserves the
canonical class (since both $\sigma$ and $-1$ act as $-1$ on it).
$X\stackrel{\pi\,}\to S$ is a $T$-fibration over $S$ with no singular
fibres, monodromy $-1$ around $\pi_1(S)$, and a section $S\into X$ at
$0\in T$. It is also a $K3$-fibration over $\Pee^1=T/\pm1$ with
four singular fibres which are double fibres modeled on $S$.

Then $\pi_1(X)$ is easily seen to be given as an extension $1\to\Z^2
\to\pi_1(X)\to\Z_2\to1$ (where the $\Z^2$ is $\pi_1(K3\times T)$). Fixing
generators $\alpha,\beta$ of $\pi_1(T)\into\pi_1(X)$ ($T\into X$ as the
fibre of $X\to\Pee^1$), and letting $\gamma\in\pi_1(S)\into\pi_1(X)$ be
the generator of $\Z_2$ ($S\into X$ as the section of $X\stackrel{\pi\,}
\to S$), we have a presentation
$$
\pi_1(X)=\langle\alpha,\beta,\gamma\rangle/(\alpha\beta=\beta\alpha,\,
\gamma^{-1}\alpha\gamma=\alpha^{-1},\,\gamma^{-1}\beta\gamma=\beta^{-1}
,\,\gamma^2=1)
$$
whose abelianisation is $H_1(X;\Z)=\Z_2^3$ (generated by $\alpha,
\beta,\gamma$). The corresponding flat line bundles, given by the
representation \{generator\,$\mapsto(-1)\in U(1)\}$ will also be
denoted by $\alpha,\beta,\gamma$. Note that $\gamma=\pi^*K_S$ is the
pull-back of the canonical bundle of $S$ to $X$.

As before pick an acyclic line bundle $L$ of degree $-1$ on $S$:
$$
H^i(L)=0=H^i(L^*)=H^i(L\otimes\gamma)\quad\forall i,
\qquad c_1(L)^2=-2,\ c_1(L)\,.\,\omega=-1,
$$
with respect to an \emph{integral} K\"ahler form $\omega$ on $S$.
(For instance on the Enriques surface studied in (\cite{BPV} V.23)
with corresponding $K3$ the double cover of $\Pee^1\times\Pee^1$ branched
over a certain $(4,4)$-curve, the line bundle $\OO(-1,1)$ restricted to
$K3$ descends to such a line bundle on $S$. The two Fubini-Study forms
on the $\Pee^1$\,s pull-back, restrict and descend to integral forms
$\omega_1,\,\omega_2$ on $S$; we then choose $\omega=\omega_1+2\omega_2$.)

Next define bundles $A$ and $B$ on $S$ by the Serre construction as before,
\begin{equation} \label{a*}
0\to\OO\stackrel{s_{A^*}\,}\To A^*\to\I_1\otimes L^*\to0,
\end{equation}
and
\begin{equation} \label{b}
0\to\OO\stackrel{s_B\,}\To B\to\I_2\otimes L^*\to0,
\end{equation}
with $\I_1$ and $\I_2$ the ideal sheaves of a point and a zero-dimensional
subscheme of length 2, respectively, where \emph{we take the point to lie
in the length 2 subscheme} this time.

We are ready to define
\begin{equation} \label{e}
A':=\pi^*A\otimes\OO(3S)\otimes\gamma, \quad B':=\pi^*B\otimes\OO(-3S)
\otimes\gamma, \quad E:=A'\oplus B',
\end{equation}
where $S\into X$ is the zero-section of $X\stackrel{\pi\,}\to S$, defining
a divisor with corresponding line bundle $\OO(S)$.

\begin{Theorem}
The rank 4 bundle $E$ defined in (\ref{e}) satisfies all but the last of
the conditions mentioned at the start of this paper, with respect to any
K\"ahler form $\Omega$ in the class of $12[S]+\pi^*\omega$.
\end{Theorem}

\begin{Proof}
By construction $\Lambda^4E\cong\OO_X$, and by Riemann-Roch, or
by lifting to $K3\times T$ and dividing by 2,
$$
\chi(E)=\chi(A')+\chi(B')=3\chi_S(A)+(-3)\chi_S(B)=3(c_2(B)-c_2(A))=3(2-1)=3,
$$
the second condition.

Since by choice of $\Omega$ both $A'$ and $B'$ have degree zero, to
prove slope-polystability it is enough to check that $A$ and $B$ are
slope-stable on $S$ with respect to $\omega$ (which was chosen to be
integral and such that $c_1(L)\,.\,\omega=-1$, remember).

Letting $P$ be one of $A^*$ or $B$, of degree 1, to check
stability we need only show that $P\otimes\eta$ has no sections for any
line bundle $\eta$ of degree less than or equal to $-1$ (this is where
the integrality of $\omega$ is used). But we have a sequence (\ref{a*},
\ref{b})
$$
0\to\eta\to P\otimes\eta\to L^*\otimes\eta\otimes\I\to0,
$$
for $\I$ some non-trivial ideal sheaf. $\eta$ has no sections since it
has degree $\le-1$, and $L^*\otimes\eta$ has degree $\le0$ so has no
sections with zeros. Since $\I$ is non-trivial this shows that $P$ has
no sections, as required.

We now turn to the $\Lambda^2E$ pairing, which (in the untwisted case)
splits as
$$
\spreaddiagramrows{-8mm}
\diagram
H^1(\Lambda^2A')\otimes H^1(B')\otimes H^1(B')\hspace{1cm} \ddrto \\
\oplus \\
\hspace{-1cm}(H^1(A'\otimes B')\otimes H^1(A')\otimes H^1(B'))^{\oplus2}
\rto & \C, \\ \oplus \\
H^1(\Lambda^2B')\otimes H^1(A')\otimes H^1(A')\hspace{1cm} \uurto
\enddiagram
$$
and similarly for the dual pairing. We compute these using the Leray
spectral sequence for $X\stackrel{\pi\,}\to S$, noting that
$\pi_*\OO(nS)\cong\OO_S^{\oplus n}\ (n\ge0)$ and so, by relative Serre
duality, $R^1\pi_*\OO(-nS)\cong\gamma^{\oplus n}$. (Again we are
suppressing some pull-backs and identifying $\gamma$ with $K_S$.)

Thus $H^i(\Lambda^2A')\cong H^i_S(L^{\oplus6})=0$ and
$H^i(\Lambda^2B')\cong H^{i-1}_S((L^*\otimes\gamma)^{\oplus6})=0$,
and the same holds on twisting by flat line bundles or taking duals.
Therefore only the central pairing above survives.

Tensoring (\ref{b}) by $A$ or $A\otimes\gamma$, and using the fact that
$H^0(A^*\otimes\I_2)=0$ (since $A^*$ has only one section $s_{A^*}$,
and this vanishes at one point only), we see that
$$
H^0(A\otimes B)=0=H^0(A\otimes B\otimes\gamma),
$$
and the same also holds with $A,B$ replaced by $A^*,B^*$. Thus the
$H^0_S(R^1\pi_*\OO_X\otimes\,.\,)$ terms that appear in the pairing (from
the Leray spectral sequence) vanish, and we are left with two copies of
$$
H^1_S(\pi_*\OO_X\otimes A\otimes B)\otimes H^1_S
((A\otimes\gamma)^{\oplus3})\otimes H^0_S(B^{\oplus 3})\to\C,
$$
and the dual pairing is twice
$$
H^1_S(\pi_*\OO_X\otimes A^*\otimes B^*)\otimes H^0_S
((A^*)^{\oplus3})\otimes H^1_S((B^*\otimes\gamma)^{\oplus 3})\to\C.
$$
Twisting by any flat line bundle that is
non-trivial on the $T$ fibres destroys the $\pi_*\OO_X$ term
(this was the original idea for the whole construction of course),
so we need only consider twisting the first terms by $\gamma$.
Since the tensor product of the three line bundles we tensor
by must be trivial, one of the other two terms must be twisted by
something containing a $\gamma$ factor (i.e. not in the span of
$\alpha,\beta$).

Thus the
pairings become, by Serre duality, 2.3.3=18 copies of
\begin{equation} \label{1}
H^1_S(A\otimes\gamma)\otimes H^0_S(B)\to H^1_S(A\otimes B\otimes\gamma),
\end{equation}
and the dual
\begin{equation} \label{2}
H^0_S(A^*)\otimes H^1_S(B^*\otimes\gamma)\to H^1_S(A^*\otimes B^*
\otimes\gamma).
\end{equation}
Twisting either of the $H^0$ groups above by $\gamma$ destroys them
by inspection of (\ref{a*},\ \ref{b}). So from what we have already
proved it is now enough to show that the pairings are non-trivial
as they stand but trivial when the $H^1$ groups above are twisted by
$\gamma$.

For the first case (\ref{1}) note that the required vanishing
follows from the vanishing of $H^1_S(A)$: $c_2(A)$ was chosen to be 1
making $\chi(A)=0=h^0(A)-h^1(A)+h^0(A^*\otimes\gamma)=-h^1(A)$. The
non-triviality of the untwisted pairing follows by taking the
cohomology of $A\otimes\gamma\otimes\,$(\ref{b}):
$$
0\to H^1(A\otimes\gamma)\stackrel{s_B\,}\To H^1(A\otimes B\otimes\gamma),
$$
where the first zero follows from $H^0(A^*\otimes\gamma)=0$ (\ref{a*}).
Since the first group has dimension 1 (also by Riemann-Roch, since
$h^2(A\otimes\gamma)=h^0(A^*)=1$) the pairing (\ref{1}) is non-zero.

For the dual pairing (\ref{2}) we take the cohomology of (\ref{a*})$\,
\otimes B^*\otimes\gamma$, giving
$$
0\to H^0(B\otimes\I_1\otimes\gamma)\to H^1(B^*\otimes\gamma)
\stackrel{s_{A^*}\,}\To H^1(A^*\otimes B^*\otimes\gamma).
$$
The pairing is the second map so we want the first map to \emph{not}
be onto, but to be onto when the $\gamma$\,s are removed (so that the
twisted pairing vanishes). Recalling that we chose the zeros of $s_{A^*}$
to lie in those of $s_B$ we see that $H^0(B\otimes\I_1)=\C\,.\,s_B$,
while $H^0(B\otimes\I_1\otimes\gamma)=0$. Since, by Riemann-Roch,
$h^1(B^*\otimes\gamma)=2$ and $h^1(B^*)=1$, this gives the required result.
\end{Proof}

Again finding a stable $F$ with $c_2(F)=c_2(X)-c_2(E)$ and $c_1(F)=0$
has defeated me. For $F$ rank 2 this works out as $c_2(F(3S)\otimes L^*)
=5[T]$, for what it's worth.

\section{$(K3\times T^2)/\Z_2$ again}

We give a final example, again satisfying all but the last of Witten's
conditions. Since the example will be very similar to the last one in
all but the pairings (as $L^*$ will no longer be acyclic) we
will concentrate mostly on them.

Pick an Enriques surface with a $-2$-sphere $C$ in it,
let $K3$ be its universal cover, and again set $X=(K3\times T^2)/\Z_2$.
We will use the notation of the last section.

Letting $L=\OO(-C)$, $L$ and $L^*\otimes\gamma$ are acyclic while
$$
H^i(L^*)=\left\{\begin{array}{cl}\C & i=0 \\ \C & i=1 \\ 0 & i=2
\end{array} \right. \qquad\mathrm{and}\qquad 
H^i(L\otimes\gamma)=\left\{\begin{array}{cl}0 & i=0 \\ \C & i=1 \\ \C
& i=2. \end{array}\right.
$$

We can then construct $A$ and $B$ on $S$ by the Serre construction, using
the fact that $H^2(L)=0$ so there is no obstruction to finding locally
free sheaves with presentations
\begin{equation} \label{a2}
0\to\OO\stackrel{s_{A^*}\,}\To A^*\to L^*\otimes\I_1\to0,
\end{equation}
and
\begin{equation} \label{b2}
0\to\OO\stackrel{s_{B}\,}\To B\to L^*\otimes\I_2\to0,
\end{equation}
with $\I_1,\,\I_2$ the ideal sheaves of 1 and 2 points in $C\subset S$,
respectively. We take the one point to be one of the two points,
as in the last example.

Setting
$$
A':=\pi^*A\otimes\OO(3S)\otimes\gamma, \quad B':=\pi^*B\otimes\OO(-3S)
\otimes\gamma, \quad E:=A'\oplus B',
$$
we can, as in the last two examples, choose compatible K\"ahler forms on
$S$ and $X$ such that $A$ and $B$ are slope-stable and $A',\,B'$ have degree
zero, so that $E$ is slope-polystable. Of course $\Lambda^4E$ is
trivial and $\chi(E)=3$.

Then
$$
H^1_X(\Lambda^2A')=H^1_X((\Lambda^2B')^*)=H^1_S(L\otimes\pi_*\OO(6S))
=H^1_S(L)^{\oplus6}=0,
$$ and
$$
H^1_X(\Lambda^2B')=H^1_X((\Lambda^2A')^*)=H^0_S(L^*\otimes
R^1\pi_*\OO(-6S))=H^0_S(L^*\otimes\gamma)^{\oplus6}=0,
$$
so that as before the pairing on $\Lambda^2 E$ reduces to its $H^1(A'
\otimes B')$ summand, and similarly for the dual pairing and the twists
by flat line bundles.

The usual arguments give the vanishing of $H^0(A\otimes B)$, its
dual and their twists by $\gamma$, so that the pairing reduces to
18 copies of
$$
H^1_S(\pi_*\OO_X\otimes A\otimes B)\otimes H^1_S(A\otimes\gamma)
\otimes H^0_S(B)\to\C,
$$
and the dual to 18 copies of
$$
H^1_S(\pi_*\OO_X\otimes A^*\otimes B^*)\otimes H^0_S(A^*)
\otimes H^1_S(B^*\otimes\gamma)\to\C.
$$
On twisting $E$ by any flat line bundle with a non-zero $\alpha$ or
$\beta$ component the $\pi_*\OO_X$ term vanishes so the pairing does
too. Thus we need only consider twists by $\gamma$.

Thus by Serre duality the pairing is equivalent to
$$
\diagram
H^1(A\otimes\gamma)^{\oplus2}\rto^{s_B\oplus s_B'\ } & 
H^1(A\otimes B\otimes\gamma),
\enddiagram
$$
where we have noted from (\ref{b2}) that $H^0(B)\cong\C^2$ and we have
picked a basis $s_B,\,s'_B$. Similarly the dual pairing is represented by
$$
H^1(B^*\otimes\gamma)\stackrel{s_{A^*}\,}\To H^1(A^*\otimes
B^*\otimes\gamma).
$$
We want these to be non-zero, but zero on twisting by $\gamma$.

For the first pairing we tensor (\ref{b2}) by $A(\otimes\gamma)$, giving
$$
0\to H^1(A(\otimes\gamma))\stackrel{s\,}\To H^1(A\otimes B(\otimes\gamma)),
$$
where $s$ is any section of $B$ (which we see from (\ref{b2}) vanishes
on 2 points in $C\subset S$). Thus it is sufficient to show that
$H^1(A\otimes\gamma)\ne0$ and $H^1(A)=0$. But $h^1(A\otimes\gamma)=
-\chi(A\otimes\gamma)+h^0(A\otimes\gamma)+h^0(A^*)=0+0+1=1$, and
$h^1(A)=-\chi(A)+h^0(A)+h^0(A^*\otimes\gamma)=0+0+0=0$.

For the dual pairing, tensoring (\ref{a2}) with $B^*(\otimes\gamma)$
yields
$$
0\to H^1(B\otimes\I_1(\otimes\gamma))\to H^1(B^*(\otimes\gamma))
\stackrel{s_{A^*}\,}\To H^1(A^*\otimes B^*(\otimes\gamma)).
$$
With the $\gamma$ the last map is the dual pairing, and by Riemann-Roch
and $h^0(B)=2,\,\chi(B)=-1$ we see the sequence is
$$
0\to0\to\C^3\stackrel{s_{A^*}\,}\To H^1(A^*\otimes B^*\otimes\gamma)
$$
so that the pairing is indeed non-zero. Removing the $\gamma$\,s
gives the twisted dual pairing, which we would like to show vanishes.

Since the zero of $s_{A^*}$ was chosen to lie in the two zeros of $s_B$,
we see that $h^0(B\otimes\I_1)\ge1$, while by Riemann-Roch $h^1(B^*)=
1+h^0(B^*)+h^0(B\otimes\gamma)=1$, so that the sequence becomes
$$
0\to\C\to\C\stackrel{s_{A^*}\,}\To H^1(A^*\otimes B^*)
$$
and the final map must be zero.

\begin{flushleft}
Department of Mathematics, Harvard University, One Oxford Street,
Cambridge MA 02138, USA. \\
\small Email: \tt thomas@maths.ox.ac.uk
\end{flushleft}

\end{document}